\documentclass{article}

\usepackage{float}

\usepackage{arxiv}
\usepackage{xcolor}
\usepackage[utf8]{inputenc} 
\usepackage[T1]{fontenc}    
\usepackage{hyperref}       
\usepackage{url}            
\usepackage{booktabs}       
\usepackage{amsfonts}       
\usepackage{nicefrac}       
\usepackage{microtype}      
\usepackage{lipsum}
\usepackage{graphicx}
\usepackage{amsmath}
\usepackage{cite}
\graphicspath{ {./images/} }
\usepackage[english]{babel}

\usepackage{stackengine}

\title{Effect of Choice of Boundary Condition on the Numerical Efficiency of Direct Solutions of Finite Difference frequency Domain Systems with Perfectly Matched Layers}

\author{
 Nathan Zhao \\
  Department of Applied Physics\\
  Stanford University\\
  Stanford, CA 94305 \\
  \texttt{nzz2102@stanford.edu} \\
  \And
 Shanhui Fan \\
  Department of Electrical Engineering\\
  Stanford University\\
  Stanford, CA 94305 \\
  \texttt{shanhui@stanford.edu} \\
}

\begin{document}
\maketitle
\begin{abstract}
    Direct solvers are a common method for solving finite difference frequency domain (FDFD) systems that arise in numerical solutions of Maxwell's equations. In a direct solver, one factorizes the system matrix. Since the system matrix is typically very sparse, the fill-in of these factors is the single most important computational consideration in terms of time complexity and memory considerations. As a result, it is of great interest to determine ways in which this fill-in can be systematically reduced. In this paper, we show that in the context of commonly used perfectly matched boundary layer methods, the choice of boundary condition behind the perfectly matched boundary layer can be exploited to reduce fill-in incurred during the factorization, leading to significant gains of up to 40\% in the efficiency of the factorization procedure. We illustrate our findings by solving linear systems and eigenvalue problems associated with the FDFD method.
\end{abstract}


\section{Introduction}
In any finite difference simulation, the domain must be truncated and a boundary condition is imposed. Common
boundaries include the periodic boundary condition or the Dirichlet boundary condition. However, in many applications, one wants to simulate an "infinite" domain in the sense that any wave incident on the boundary should have no reflection. In one dimension, it is possible to eliminate all reflections via the Mur
absorbing boundary condition [1]. However in $n > 2$ dimensions, there is no absorbing boundary condition that can eliminate reflections for waves incident from all angles. Instead, it is common to surround the computational domain by perfectly matched layers (PML), a specially designed anisotropic absorber \cite{Berenger1994, taflove2005}. The PML region is typically truncated using a boundary condition such as the periodic or the Dirichlet boundary condition. 

The PML has been extensively used in computational electromagnetics, for both frequency and time-domain solutions of the Maxwell’s equations \cite{Berenger1994, Berenger1996, sacks1995, johnson2021, taflove2005}. In the time domain, significant efforts have been invested in reducing the reflection for a given number of layers in the PML region \cite{Collino1998,Agrawal:04, oskoii2011}. Various formalisms of the PML have also been developed to ensure strong suppression of reflection over broad bandwidth \cite{Berenger1996_fdtd, taflove2005}. In the frequency domain, there are additional considerations in the applications of PML. For example, Ref. \citenum{Wonseok2012} has shown that the condition number of the system matrix, and hence the performance of iterative solutions of Maxwell’s equations in the frequency domain, can vary drastically for different formalisms of PML, even though these formalisms offer similar wall-clock performance in time-domain algorithms. 

In this paper we provide a discussion of PML in the context of direct solver for Maxwell’s equations in the frequency domain. In such direct solvers, one typically performs an LU decomposition of the system matrix. We show that the fill-in in the LU decomposition strongly depends on the boundary condition that is used to truncate the PML. In particular, the use of the Dirichlet boundary condition leads to far less fill-in and hence significant efficiency gain, as compared with the use of the periodic boundary condition. This result is in contrast with both time-domain algorithms, as well as with iterative algorithms in the frequency domain for the solutions of Maxwell’s equations, where the use of these two boundary conditions do not significantly affect the performance of the algorithms. Our result is important for the development of direct solvers for Maxwell’s equations in the frequency domain.

\section{A Brief Review of the Finite Difference Frequency Domain Method}
\subsection{Linear system formulation}
The Maxwell's equation written for the $\mathbf{H}$-field is: 
\begin{equation}
\nabla \times \frac{1}{ \overleftrightarrow{\epsilon_r}(\mathbf{r})} \nabla \times \mathbf{H(r)}   - \omega^2 \mu_0 \mathbf{H(r)} = -i\omega \mathbf{M(r)}
\label{eq:maxwell}
\end{equation}
where $\overleftrightarrow{\epsilon_r}$ denotes a three-dimensional dielectric tensor, $\omega$ is the frequency, $\mathbf{H(r)}$ denotes the vector magnetic field, $\mathbf{M(r)}$ is the vector magnetic source, and the magnetic permeability is assumed to be $\mu_0$ in the entire computational domain. In the finite-difference frequency-domain (FDFD) method, one discretizes Eq. \eqref{eq:maxwell} via finite differences on a Yee grid \cite{yee1966, Veronis2004, taflove2005, leveque2007} to obtain a system of linear equations:
\begin{equation}
      \hat A \mathbf{x} =\mathbf{b} 
    \label{eq:forward}
\end{equation}
The matrix $\hat A$ corresponds to the left hand side of Eq. \eqref{eq:maxwell} and $\mathbf{x}$ is a vector corresponding to the discretized version of $\mathbf{H}$. $\mathbf{b}$ corresponds to the source term on the right hand side.

In order to obtain a solution to Eq. \eqref{eq:forward}, a number of algorithms exist, such as Gaussian elimination \cite{strang2016}. However, a common method is to factorize the matrix $\hat A$ into $\hat A = \hat L \hat U$,  where $\hat L$ is a lower triangular matrix and $\hat U$ is an upper triangular matrix. The solution can then be obtained:
\begin{equation}
    \mathbf{x} = \hat U^{-1}(\hat L^{-1} \mathbf{b})
    \label{eq:lu_soln}
\end{equation}
Eq. \eqref{eq:lu_soln} is typically called back-substitution. The computation of $L^{-1} b$ and $U^{-1} (L^{-1}b)$ is cheaper as compared with the full Gaussian elimination of the system matrix $\hat A$  \cite{Li2005,Davis2004}. However, the process of determining the factorization requires a modified version of Gaussian elimination, so the overall time complexity of using the LU method is the same as the Gaussian elimination \cite{strang2016}. On the other hand, if one has multiple right hand sides to solve, then the LU factorization becomes efficient as the $L$ and $U$ factors can be re-used.

\subsection{Linear Eigenvalue System Formulation}\label{sec:linear_eigen}

The LU decomposition is also useful for solving eigenvalue problems in electromagnetics, an example of which is:
\begin{equation}
\nabla \times \frac{1}{ \overleftrightarrow{\epsilon_r}(\mathbf{r})} \nabla \times \mathbf{H(r)}   = \omega^2 \mu_0 \mathbf{H(r)}
\label{eq:maxwell_eigen}
\end{equation}
Following the same discretization procedure as that leading to Eq. \eqref{eq:forward}, we obtain the eigenvalue problem:
\begin{equation}
\hat A \mathbf{x} = \lambda \mathbf{x}
    \label{eq:eigen}    
\end{equation}
where $\hat A$ now corresponds to the left-hand side  of Eq. \eqref{eq:eigen} and $\lambda = \mu_0\omega^2$. Eigenvalues are typically determined iteratively by using the Arnoldi iteration and its variants where one repeatedly applies the matrix $\hat A$ to a test vector \cite{Lehoucq97arpackusers, golub2000}. However, direct application of these methods can only produce the largest $k$ eigenvalues of a system. In computational electromagnetics, one typically is interested instead in eigenvalues that are near a target value $\sigma$. For such a problem, a reformulation of the eigenvalue problem is required \cite{Lehoucq97arpackusers,Campos2012}. Instead of solving Eq. \eqref{eq:eigen} we solve:
\begin{equation}
    (\hat A-\sigma \hat I )^{-1}\mathbf{x} = \beta \mathbf{x} 
    \label{eq:shift_invert_mode}
\end{equation}
where $\hat I$ represents the identity matrix. The eigenvalues $\beta$ can be related to eigenvalues of $A$ via: $\beta = 1/(\lambda - \sigma)$. Thus the eigenvalues of $\hat A$ that are closest to $\sigma$ correspond to the largest eigenvalues of $(\hat A-\sigma I)^{-1}$. In applying iterative methods to obtain the largest eigenvalues of $(\hat A-\sigma I)^{-1}$, we need to repeatedly apply $(\hat A-\sigma \hat I)^{-1}$. In the FDFD method, $\hat A$ is sparse, but $(\hat A- \sigma \hat I)^{-1}$ is fully dense. For large problems it is too costly memory-wise to form $(\hat A-\sigma \hat I)^{-1}$ explicitly. The efficient way to do this is to compute an LU factorization of $\hat A-\sigma \hat I$ \cite{Peng1996}. Thus, the LU factorization is important for the eigenvalue problems that arise in the FDFD method as well. 

\subsection{Brief Review of Boundary Conditions}\label{sec:boundary_intro}

\begin{figure}[H]
    \centering
    \includegraphics[width = 5 in]{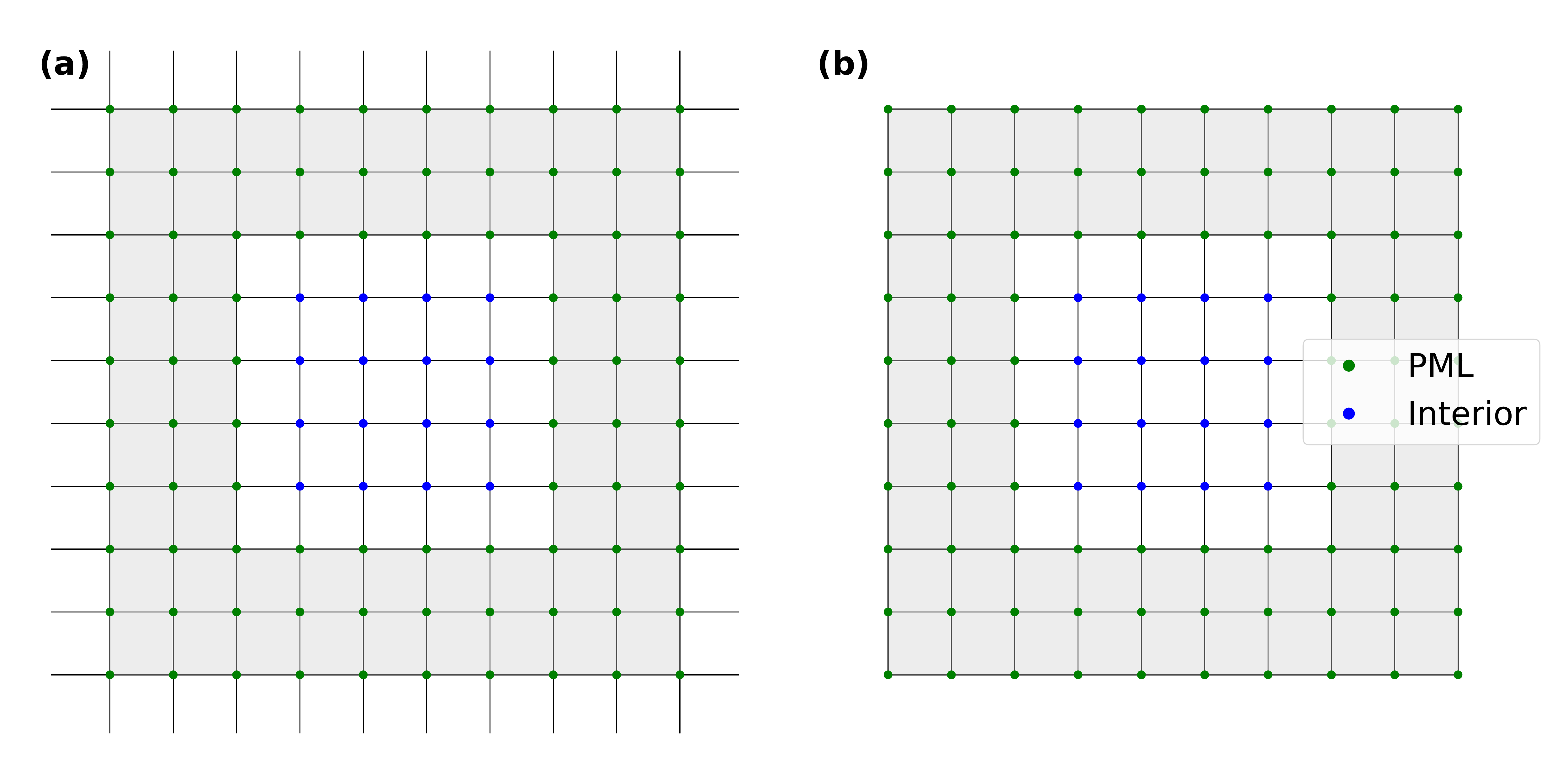}
    \caption{Two dimensional finite difference grid with a five point stencil. a) A grid truncated with the periodic boundary condition. The boundary connections which connect left boundary to right and top to bottom are shown as trailing black edges. b) Same grid with the Dirichlet boundary condition.}
    \label{fig:pbc_to_pec}
\end{figure}

In the FDFD method the computational domain is finite and a boundary condition needs to be applied on the surface of the computational domain. This is illustrated in  Fig. \ref{fig:pbc_to_pec}, where we show an undirected graph representation of a simple square grid consisting of $N = N_x \times N_y$ nodes with $ N_x = N_y = 10$. The points at which the fields are evaluated are termed nodes. We denote the fields evaluated at the nodes as $\mathbf{x}(i,j)$, where $i$ is the row index and $j$ is the column index. For a 2D problem, the Maxwell's equations can be written in terms of only $E_z$ or $H_z$ fields, where $z$ is the out-of-plane direction, hence $\mathbf{x}$ can correspond either to the $E_z$ or the $H_z$ fields. The outermost nodes are the boundary nodes where the boundary condition is enforced. The edges represent the node-to-node coupling as derived from finite-difference approximation of Eq. (1). In Fig. \ref{fig:pbc_to_pec}a, we show the periodic boundary, which can be described as: $\mathbf{x}(N_x,j)=\mathbf{x}(0,j)$ and $\mathbf{x}(i,N_y) = \mathbf{x}(i,0)$. In Fig. \ref{fig:pbc_to_pec}b, we show a Dirichlet boundary condition:  $\mathbf{x}(N_x,j)=\mathbf{x}(0,j) = \mathbf{x}(i,N_y) = \mathbf{x}(i,0) = 0$. When $\mathbf{x}$ corresponds to the $\mathbf{H}$($\mathbf{E}$) fields, the Dirichlet boundary condition corresponds to the perfectly magnetic (electric) conductor boundary condition. 

In typical FDFD simulations one often surrounds a computation domain of interest with perfectly match layers (PML) \cite{Berenger1994, Wonseok2012}. These layers are used to ensure that any wave exiting the computational domain is nearly perfectly absorbed, which is important for simulating an open electromagnetic system. When the PML are used, the entire computational domain consists of the domain of interest and the PML regions. Outside the PML regions, the computational domain still needs to be truncated. Either the periodic or the Dirichlet boundary conditions can be used for this purpose. 
 
The choice of either the periodic or the Dirichlet boundary conditions outside the PML regions does not affect their absorption properties. Moreover, in the finite difference time domain method, or when iterative solvers are used in the finite difference frequency domain method, these two boundary conditions result in very similar computational cost. Therefore, there has been very little attention paid to the choice between these two boundary conditions. As we will show in this paper, however, this choice in fact makes a significant difference when direct solvers are used in the finite difference frequency domain method, since the choice significantly influences the fill-in behaviors in the LU decomposition of the system matrix.  



\section{A Brief Review of the Fill-in of Sparse Linear Systems}\label{intro_fill}\label{sec:heuristics}
In general, $\hat A$ for Eq. \eqref{eq:forward} or Eq. \eqref{eq:eigen} is sparse. A key consideration in the LU factorization is the issue of fill-in, which is the accumulation of additional non-zero elements in $\hat L$ and $\hat U$ outside those of $\hat A$. In general, fill-in during the factorization process for a sparse linear system depends sensitively on the detailed sparsity pattern of the system matrix.  Consider, for example, the following factorization of a matrix:
\begin{equation}
\begin{bmatrix}
    p & p& p & p & p  \\
    p & p&   &   &   \\
    p &  & p &   &   \\
    p &  &   & p &  \\
    p &  &   &   & p \\
    \end{bmatrix} = \begin{bmatrix}
     p &  &  &  &   \\
     p & p &  &  &   \\
     p & p & p & &   \\
     p & p & p & p&  \\
     p & p & p & p & p \\
    \end{bmatrix} \begin{bmatrix}
    1 & p& p &  p& p   \\
      & 1& p & p & p  \\
      & &  1 & p& p  \\
      &  &  &  1& p \\
      &  &  &   & 1 \\
    \end{bmatrix}
    \label{eq:suboptimal_ordering}
\end{equation}
where $p$ denotes a matrix element which is nonzero (element values may be different). Because the first row has no zeros, the factorization process results in $\hat L$ and $\hat U$ factors that are dense \cite{strang2016}. In contrast, consider the factorization of an equivalent system as obtained from the left hand size of Eq. \eqref{eq:suboptimal_ordering} by permuting the row index:
\begin{equation}
    \begin{bmatrix}
    p & &  &  & p  \\
      & p&  &  & p  \\
      & & p & &  p \\
      &  &  & p& p \\
     p &p  & p & p& p \\
    \end{bmatrix} =\begin{bmatrix}
    p & &  &  &    \\
      & p&  &  &   \\
      & & p & &   \\
      &  &  & p&  \\
     p &p  & p & p& p \\
    \end{bmatrix} \begin{bmatrix}
    1 & &  &  & p  \\
      & 1&  &  & p  \\
      & & 1 & &  p \\
      &  &  & 1& p \\
      &  &  &  & 1 \\
    \end{bmatrix} 
    \label{eq:optimal_ordering}
\end{equation}
The $\hat L$ and $\hat U$ factors resulting from Gaussian elimination do not have any fill-in.  In general, the order of the equations as well as the order of the variables can drastically change the fill-in. The operations which generate these reorderings can be expressed as preconditioners on $\hat A$ via a set of permutation matrices $\hat P$ and $\hat Q$: $\hat A \rightarrow \hat P\hat A\hat Q$, where $\hat P$ performs row-wise and $\hat Q$ performs column-wise permutations. Unfortunately, the problem of finding the optimal permutation to produce the least amount of fill-in is NP-complete \cite{Yannakakis1981}, which means there are no known, at-worst polynomial time algorithms to determine the optimal permutation. Instead, there is a wide array of different heuristic or approximate methods to reduce fill-in. 


The state-of-the-art re-orderings are based on the exact minimum degree algorithm and their variants, such as the methods of approximate minimum degree (AMD) and nested dissection \cite{Amestoy1996, George1973, Gilbert1986,Agarwal2012, METIS}. 
All of these methods effectively rely on the minimum degree heuristic, which is the observation that first eliminating variables in the factorization process with fewer couplings (or nonzeros in the row/column) tends to yield less fill-in overall.

\section{Reducing Fill-in by Reducing Connectivity within the PML}
As discussed in Section \ref{sec:boundary_intro}, the choice of boundary condition backing a PML can vary. As a result, we have the freedom to choose the boundary condition behind the PML. In this section, we now show how the choice of boundary condition behind the PML can affect the efficiency of sparse direct solvers, specifically impacting the fill-in during factorization. 

From Figure \ref{fig:pbc_to_pec}, we can apply the minimum degree heuristic by counting the neighbors of the individual nodes on the grid. The number of neighbors corresponds to the number of nonzero elements in each row of $\hat A$. For the periodic boundary case shown in Fig. \ref{fig:pbc_to_pec}a, all nodes including the boundary have 4 neighbors. By the minimum degree heuristic, there is no preference in picking which rows to eliminate first. Meanwhile, for the Dirichlet boundary in Fig. \ref{fig:pbc_to_pec}b, the four corner nodes have 2 neighbors, the edge nodes have 3 neighbors and then all other nodes have 4 neighbors. By the minimum degree heuristic in this case, we can see that the rows corresponding to the boundary nodes should be eliminated first. 

Finally, we can take advantage of the fact that the nodes adjacent to the Dirichlet boundary have fields that are generally close to zero, since the fields have been strongly by the PML. As a result, we can also further modify the connectivity to minimize the number of neighbors per node, as we show in Fig. \ref{fig:xtreme}:

\begin{figure}[H]
    \centering
    \includegraphics[width = 5 in]{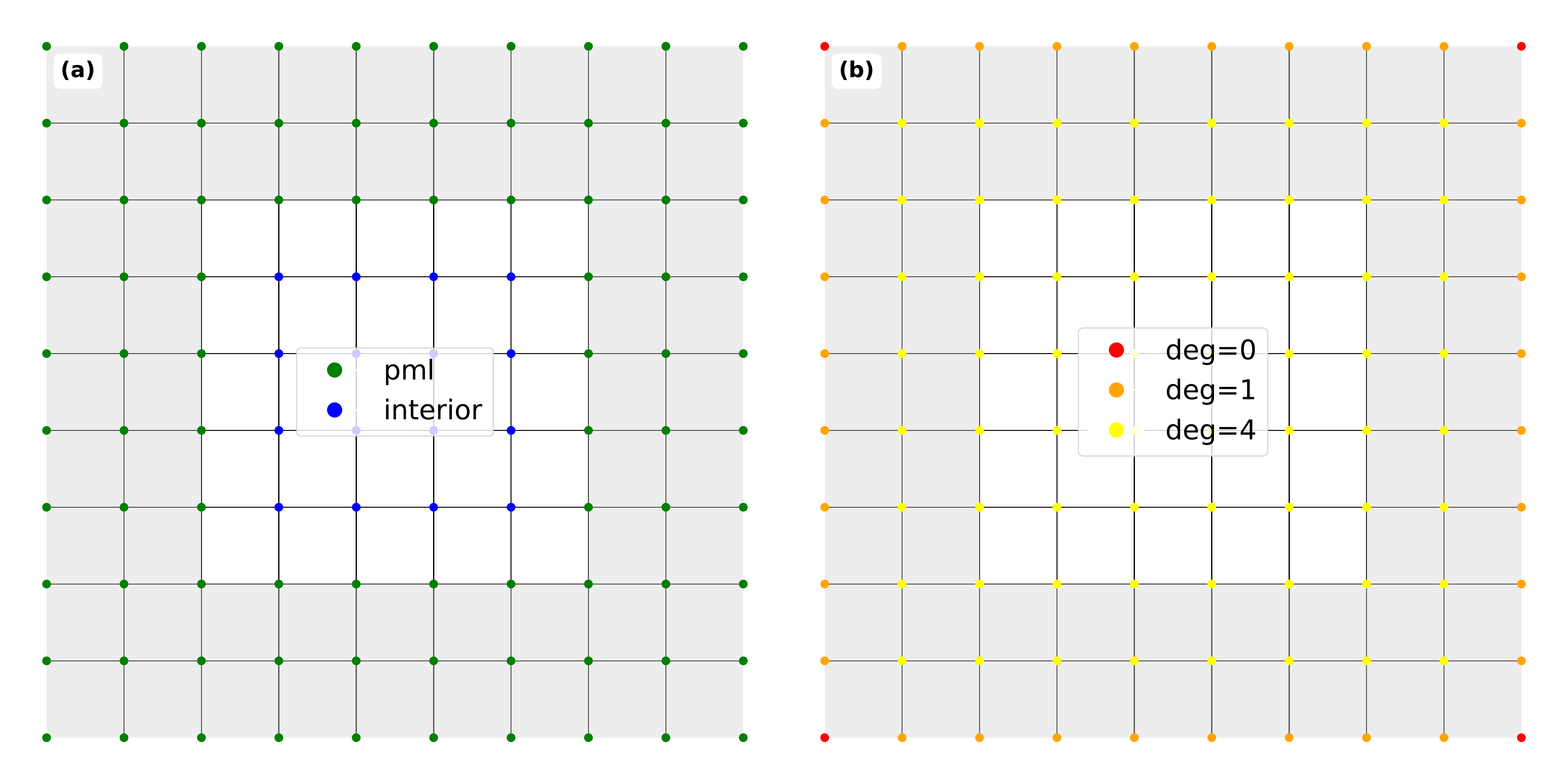}
    \caption{a) Optimized grid with reduced connectivity of the boundary. b) Same as a) with the nodes color coded by their degree (deg) in the graph representation.}
    \label{fig:xtreme}
\end{figure}

\begin{table}[H]
\begin{center}
\begin{tabular}{ |c|c|c| } 
 \hline
  Boundary Condition of Grid & Total Grid Couplings for ($N_x = N_y$) & $N_x=100$ \\ \hline 
 Periodic &  $4N_x^2$ & 40000 \\ 
 Dirichlet & $4(N_x-2)^2+4(N_x-2)+8$ & 38816 \\ 
 Modified Dirichlet & $4(N_x-2)^2+ N_x-2$ & 38514 \\ 
 \hline
\end{tabular}
\vspace{0.1 in}
\caption{Changes in the number of couplings in a grid of $N$ nodes for different boundary conditions. The 2nd column shows the analytic expression for the number of total couplings in the grid. The third column shows the values of these expressions for the specific case of $N=10000$.}
\end{center}
\label{Tab:driven_specs}
\end{table}

In Fig. 3, we show the influence of the boundary conditions on the fill-in behavior for the small 10$\times$10 grids shown in Figs. 1 and 2. The use of such a small grid allows easier visualization of various quantities involved.  We solve a 2D problem corresponding to setting all terms involving $\partial_z$ to zero in Eq.  \eqref{eq:maxwell}. We simulate a magnetic source in vacuum for the magnetic field perpendicular $H_z$ to the 2D plane:
\begin{equation}
    \bigg(\frac{\partial}{\partial x}\frac{1}{\epsilon_{r}(x,y)} \frac{\partial}{\partial x}+    \frac{\partial}{\partial y}\frac{1}{\epsilon_{r}(x,y)} \frac{\partial}{\partial y} + \mu_0 \omega^2\bigg) H_z = i\omega M_z
\label{eq:TM}
\end{equation}
where we assume $\epsilon_{r}$ is isotropic, $M_z$ is the $z$-component of the magnetic source current. We consider the periodic boundary conditions, the Dirichlet boundary conditions, and the modified Dirichlet boundary conditions as discussed above where the connectivity at the boundary is reduced. The first column of Figure 3 shows the sparsity pattern of the system matrices for the three boundary conditions. The sparsity patterns of the three matrices are similar since they differ only on the boundary. Compared with the system matrix from the Dirichlet boundary condition (Fig. 3d), the matrix from the periodic boundary condition (Fig. 3a) has extra non-zero elements at the lower-left and upper-half corners in the figure due to the long-range coupling induced by the periodic boundary condition. The matrix from the modified Dirichlet (Fig. 3g) boundary condition has less non-zero elements near the diagonal as results from the modifications.

The second column in Figure 3 shows the $\hat L$ factor as obtained by directly factoring the corresponding matrices in the first column as is without any reordering. For the case with periodic boundary condition (Fig. 3b), the $\hat L$ factor fills in near the diagonal and the lower edge in the figure. The fill-in at the lower edge arises from the long-range coupling as induced by the periodic boundary condition. By switching to the Dirichlet boundary condition (Fig. 3e), the fill-in near the lower edge disappears and the number of non-zero elements in the $\hat L$ factor significantly decreases. Modifying the Dirichlet boundary condition further reduces the fill-in and the number of non-zero elements (Fig. 3h).

The third column in Figure 3 shows the $\hat L$ factor as obtained by factoring the corresponding matrices in the first column using a minimum degree heuristics to reorder the matrix \cite{Amestoy1996}. Compared with the second column, the use of the minimum degree heuristics significantly reduces the number of fill-in for all three cases. Also, the number of non-zero elements reduces as we go from the case with the periodic boundary condition (Fig. 3c), to the case with the Dirichlet boundary condition (Fig. 3f), to the case with the modified Dirichlet boundary condition (Fig. 3i). The $\hat L$ factor shown in Fig. 3i, which was obtained for the modified Dirichlet boundary condition using the minimum degree heuristics, has 456 non-zero element, as compared with the number of non-zero elements of 500 in the system matrix with the periodic boundary condition. The result in Figure 3 provides a direct demonstration that the use the modified Dirichlet boundary condition, in combination with the minimum degree heuristics, can significantly reduce the fill-ins.

\begin{figure}[H]
    \centering
    \includegraphics[width = 4.5 in]{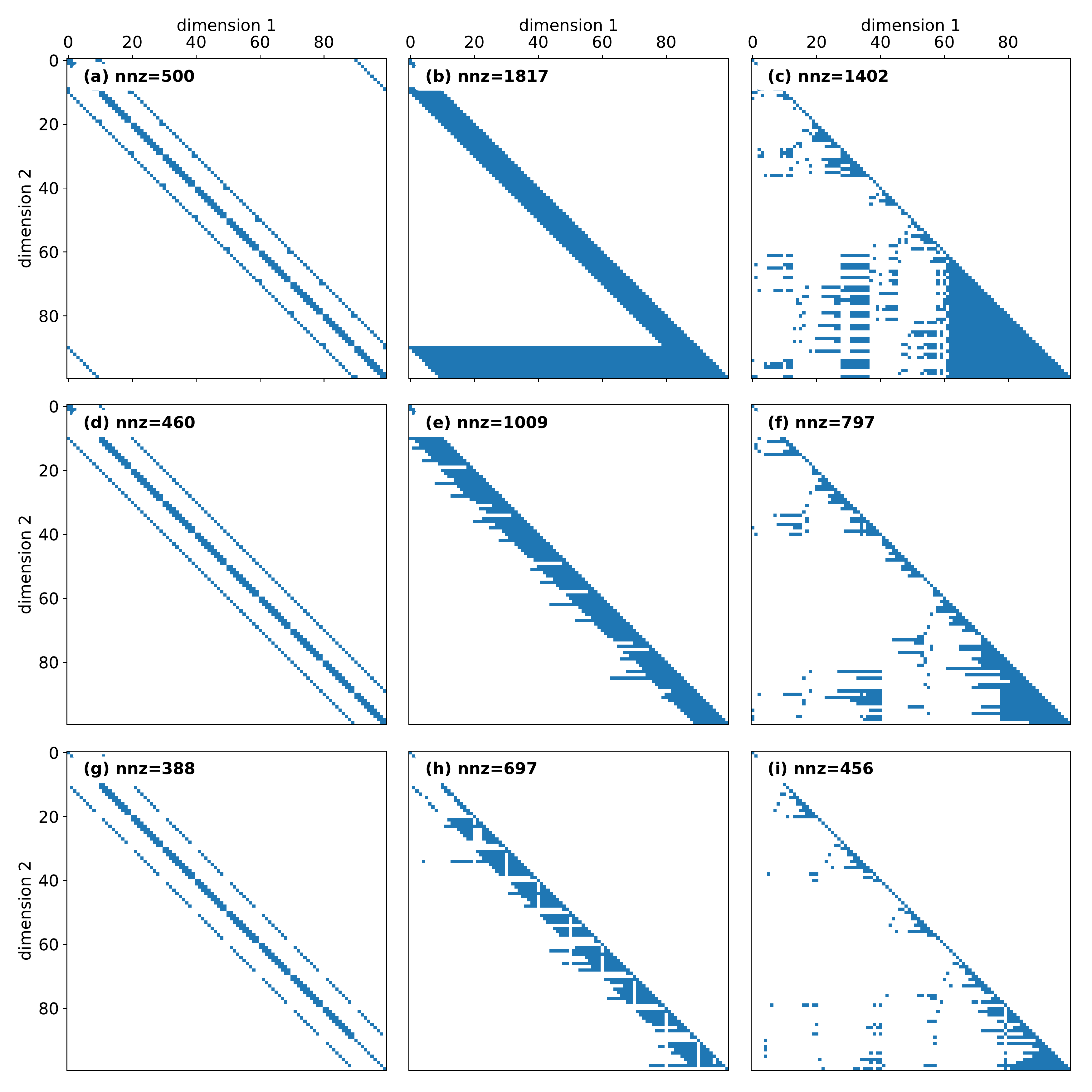}
    \caption{In the first column of Fig. \ref{fig:small_sparsity_pattern}, we show the sparsity pattern of $\hat A$ corresponding to the finite difference discretization of the $10\times 10$ grid with the a) periodic, d) Dirichlet boundary, and g) modified Dirichlet boundary with reduced connectivity at the boundary respectively. In the second column, we show the corresponding $\hat L$ factors for factorizing these three matrices as is. In the third column, we show corresponding $\hat L$ factors for factorizing these three matrices using a minimum degree heuristic to reorder the matrices. The number of nonzero elements (nnz) for each matrix is also provided in each figure.}
    \label{fig:small_sparsity_pattern}
\end{figure}

\section{Numerical Experiments}
In this section, we show that the concept as described above can be applied to typical simulations of electromagnetic waves. We consider both the driven case with a source, and eigenmode calculations without a source. For both cases, we show that choosing either the periodic or the Dirichlet boundary conditions behind the PML does not affect the accuracy of the results, but the choice of the Dirichlet boundary condition leads to better numerical performance. 

\subsection{Driven Problem}
We solve Eq. \eqref{eq:TM} as in the previous section but now the interior domain consists of $N_x \times N_y$ with $N_x = N_y = 301$ grid cells surrounded with a PML with a thickness of 30 cells to ensure minimal reflection \cite{alma2013}. The domain has a length $a = 4\lambda$ where $\lambda$ is the wavelength of the source. 

%



\begin{figure}
    \centering
    \includegraphics[width = 4.5 in]{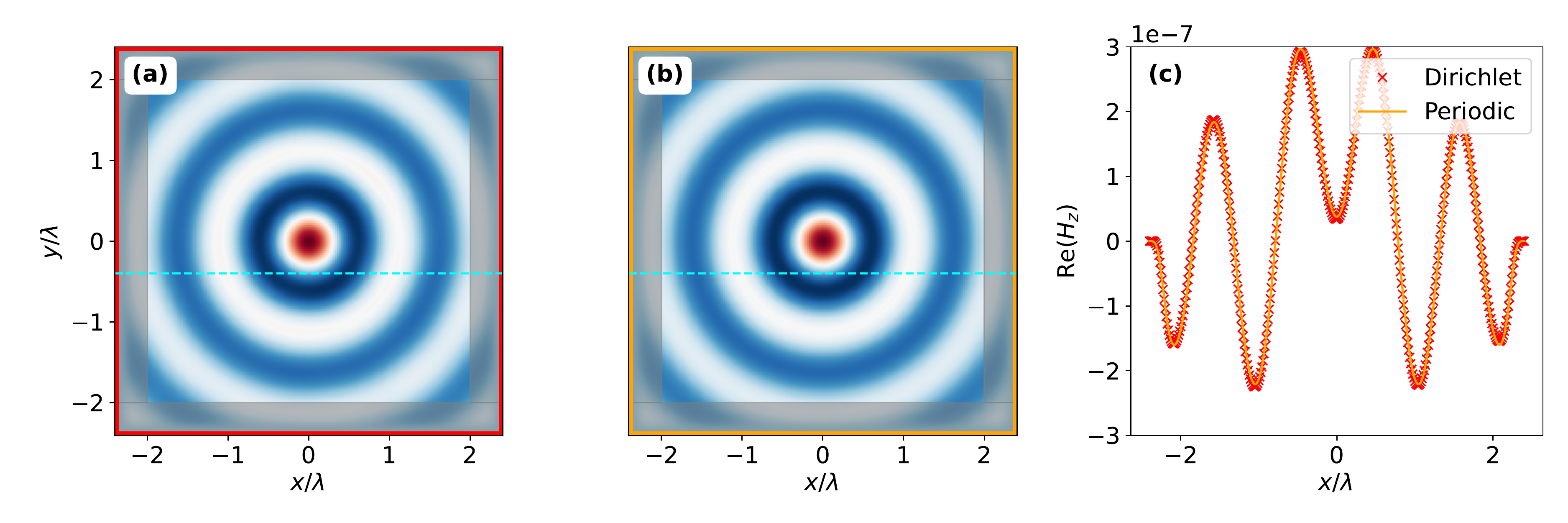}
    \caption{a) Field solution for $H_z$ for Eq. \eqref{eq:TM} with the periodic boundary. The periodic boundary is denoted by a red line. b) Field solution on the same grid but with the Dirichlet boundary condition (shown in orange). c) Comparison of the field distribution along the cyan lines in  a) and b). The transparent gray regions in (a) and (b) denote where the PML lies.}
    \label{fig:driven_validation}
\end{figure}

First, we confirm that the use of either the periodic boundary or the modified Dirichlet boundary in conjunction with a PML indeed has no effects on the solution within the interior domain. In Fig. \ref{fig:driven_validation} we show the result where the interior domain consits of vacuum and we put a point source at the center of the interior domain. For the two boundary conditions, the field distributions inside the computational domain are almost identical, as can be seen by comparing visually the field distributions shown in Fig. 4a and b, or by exampling the field distributions along a line in the computational domain as shown in Figure 4c. 

\begin{figure}
    \centering
    \includegraphics[width = 5 in]{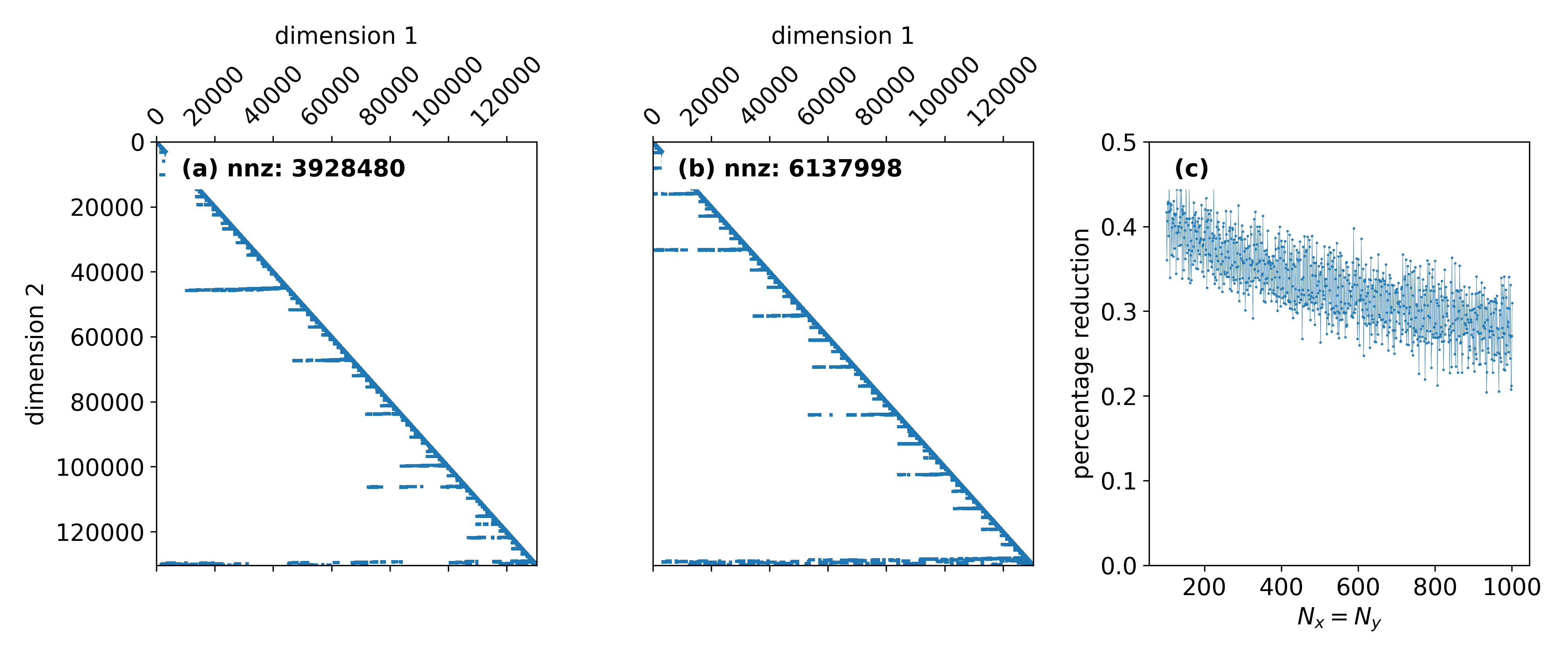}
    \caption{a) Fill-in pattern of the $\hat L$ factor for the system shown in Fig. \ref{fig:driven_validation} with the Dirichlet boundary. b) Fill-in pattern for the $\hat L$ factor for the system with the periodic boundary.  c) The percentage in the reduction of the number of nonzeros in $\hat L$ factor in going from the periodic to the modified Dirichlet boundary condition for varying values of $N_x = N_y$.}
    \label{fig:driven_sparsity_domain_scan}
\end{figure}

For the matrix $\hat A$ corresponding to the driven problem, we show in Fig. \ref{fig:driven_sparsity_domain_scan} the difference in factorizing the matrix with the periodic boundary vs the Dirichlet boundary. In Fig. \ref{fig:driven_sparsity_domain_scan}a and b, we show the fill-in of the $\hat L$ factor with the modified Dirichlet and periodic boundary respectively for the system shown in Figure \ref{fig:driven_validation}. The $\hat L$ factor from the modified Dirichlet boundary is noticeably sparser than the $\hat L$ corresponding to the periodic boundary. 

In Fig. \ref{fig:driven_sparsity_domain_scan}c, we vary the size the domain by keeping the domain to be a square shape and varying $N_x=N_y$. We plot the percentage reduction in the number of non-zero elements in the $\hat L$ factor as we change from the periodic to the modified Dirichlet boundary condition. The percentage reduction fluctuates significantly as $N_x$ varies, since the reordering strategies are heuristic. But in general, we see significant reduction for all $N_x$’s. On average the percentage reduction decreases as $N_x$ increases. This is expected since the boundary, being a smaller and smaller fraction of the domain as $N_x$ increases, would lead to a smaller effect on the relative fill-in. What is remarkable is that the boundary couplings clearly have a disproportionately large effect on the fill-in of $\hat L$ factor, even though the boundary nodes are often just 1\% or less of the total nodes in the entire grid. Even with $N_x\approx1000$, there are cases where the use of Dirichlet boundary condition results in 33\% fewer nonzeros in the $\hat L$ as compared with the use of the periodic boundary condition.

\subsection{Eigenvalue Problem}
As a second example, we consider an eigenvalue problem that commonly arises in the study of nanophotonic waveguides \cite{Veronis:05, Nicolet2006}. In this problem, we consider a structure that is periodic along the $x$-direction. For a given operating frequency $\omega$, we solve for the propagating wavevector $k_x$. To develop the formalism, we substitute via the Bloch theorem $H_z(x,y) \rightarrow H_z(x,y)e^{ik_xx}$ in Eq. \eqref{eq:TM} to obtain: 
\begin{equation}
   \frac{1}{\mu_0} \bigg( \partial_x\bigg(\frac{1}{\epsilon_{r}}\bigg)\partial_x +\partial_y\bigg(\frac{1}{\epsilon_{r}}\bigg)\partial_y +\omega^2\mu_0 +i k_x \bigg(\partial_x\epsilon_{r}^{-1}+\epsilon_{r}^{-1}\partial_x\bigg) \bigg) H_z e^{ikx} = \frac{1}{\mu_0\epsilon_{r}}k_x^2 H_z e^{ikx}
    \label{eq:maxwell_quadratic_eigen}
\end{equation}
which is a quadratic eigenvalue problem of the form:
\begin{equation}
    \hat M \mathbf{v} + \lambda^2 \hat K \mathbf{v} +\lambda \hat C \mathbf{v} = 0
    \label{eq:quadratic_eigen}
\end{equation}
where $\lambda = k_x$ and:
\begin{equation}
    \hat K = \frac{1}{\mu_0}\epsilon_r
\end{equation}
\begin{equation}
    \hat M =  \frac{1}{\mu_0}\bigg(\partial_x\bigg(\frac{1}{\epsilon_{r}}\bigg)\partial_x +\partial_y\bigg(\frac{1}{\epsilon_{y}}\bigg)\partial_y+\omega^2\bigg)
\end{equation}
\begin{equation}
    \hat C = \frac{i}{\mu_0}\bigg(\partial_x\epsilon_{r}^{-1}+\epsilon_{r}^{-1}\partial_x\bigg)
\end{equation}
Eq. \eqref{eq:quadratic_eigen} can be reformulated into an equivalent linear, generalized eigenvalue problem of the form:
\begin{equation}
    \hat D \lambda \mathbf{w} = \hat B \mathbf{w}
    \label{eq:gen_eigen_prob}
\end{equation}
where $\mathbf{w} = [\mathbf{x},\lambda \mathbf{x}]$ with $\mathbf{x}$ denoting the vector formed by $H_z(x,y)$ and where:
\begin{equation}
    \hat D = \begin{bmatrix}
        \hat K & \hat 0 \\
        \hat 0 & \hat I  \\
    \end{bmatrix}
    \label{eq:hat_D}
\end{equation}
and
\begin{equation}
    \hat B = \begin{bmatrix}
        \hat C &  \hat M \\
        -\hat I & \hat 0 \\
    \end{bmatrix}
    \label{eq:hat_b}
\end{equation}
After discretization using a computational domain consisting of a single unit cell of the waveguide structure, $\hat 0$ denotes a matrix of zeros of dimension $N_xN_y$.  In the waveguide eigenvalue problem, similar to the discussion in Section (2.2), we are typically interested in modes with wavevector near a particular value $\sigma$. Therefore, following the same derivation as in Section (2.2), the matrix that needs to factored is $\hat B-\sigma \hat D$ \cite{Lehoucq97arpackusers}.




As an illustration we consider the waveguide structure schematically shown in Figure 6a \cite{Zhao2019_meta, Zhao:20}. The permittivity of the materials used is provided in the caption of Figure 6a. The waveguide consists of two walls. Each wall consists of a mixture of dielectric and a Drude metal. The thickness of each wall is 0.3 $\mu$m. The separation between the walls is 1.6 $\mu$m.  The structure is periodic along the $x$-direction. The periodicity along the $x$-direction is 0.2 $\mu$m. The metal region in the wall has a width along the $x$-direction of 0.04 $\mu$m. Two unit cells are shown in Figure 6a. By operating this system at $\lambda = 2$ $\mu$m, the walls become near-completely reflective, creating a hollow-core waveguide with light guided in the vacuum region between the walls.

To simulate this structure we use a computational domain consisting of a single unit cell shown in Figure 6a. We discretize the computational domain with $N_x \times N_y$ grid cells where $N_x = 100$ and $N_y = 150$. Along the $x$-direction a periodic boundary condition is imposed at the edges of the computational domain. Along  the $y$-direction we place PML regions at both edges. The PML regions consist of 10 grid cells. Outside the PML region we truncate the grid with either the modified Dirichlet or the periodic boundary conditions. 

\begin{figure}[H]
    \centering
    \includegraphics[width = 5 in]{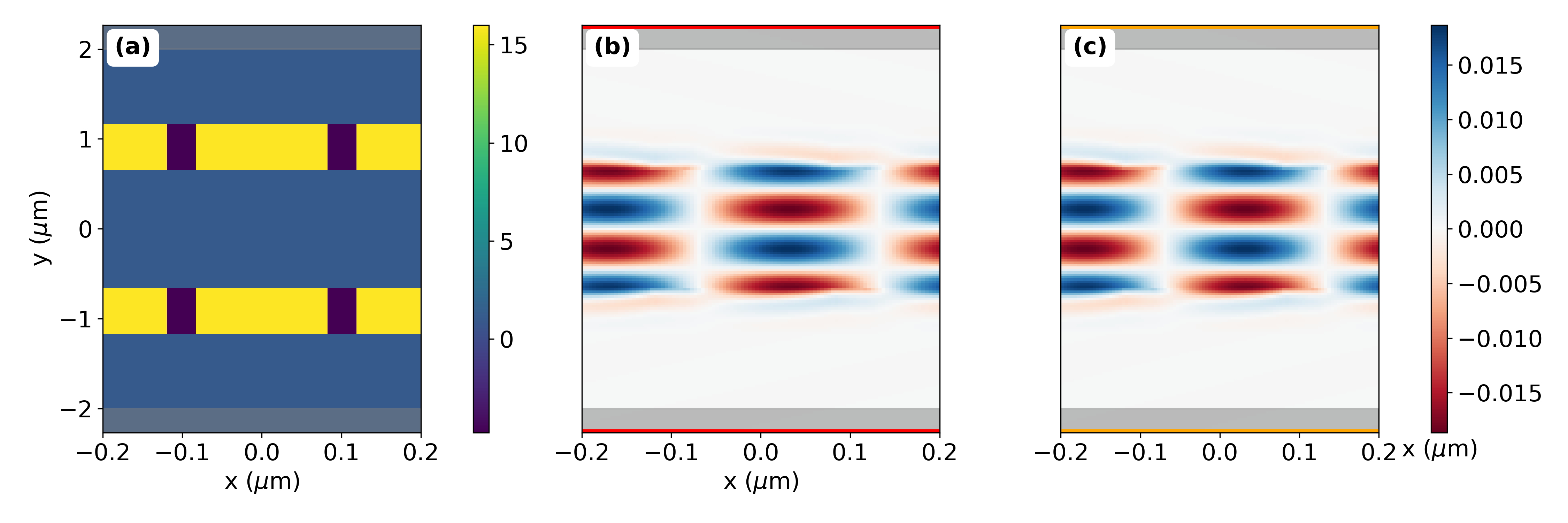}
    \caption{a) Schematic of the waveguide structure and the computational domain. The violet regions consist of a metal described by the Drude model with $\epsilon(\omega) = 1 - \omega^2/(\omega_p^2-i\gamma \omega)$, with $\omega_p = 0.72\pi \times 10^{15}$ rad/s, $\gamma = 5.5\times 10^{12}$  $s^{-1}$.  The yellow regions consist of a dielectric material with a relative permittivity of $16$. The blue region is vacuum. The gray regions are the PML regions. The structure is periodic along the $x$-direction. Plotted here are two unit cells of the structure.  b) $Re(H_z)$ modal profile of the structure with a Dirichlet boundary at the top and bottom edges of the computational domain at an operating frequency corresponding to a free space wavelength of 2 $\mu$m.  b) The corresponding result with the periodic boundary at the top and bottom edges of the computational domain.}
    \label{fig:mode_profiles_check}
\end{figure}
\begin{figure}
    \centering
    \includegraphics[width = 4.5 in]{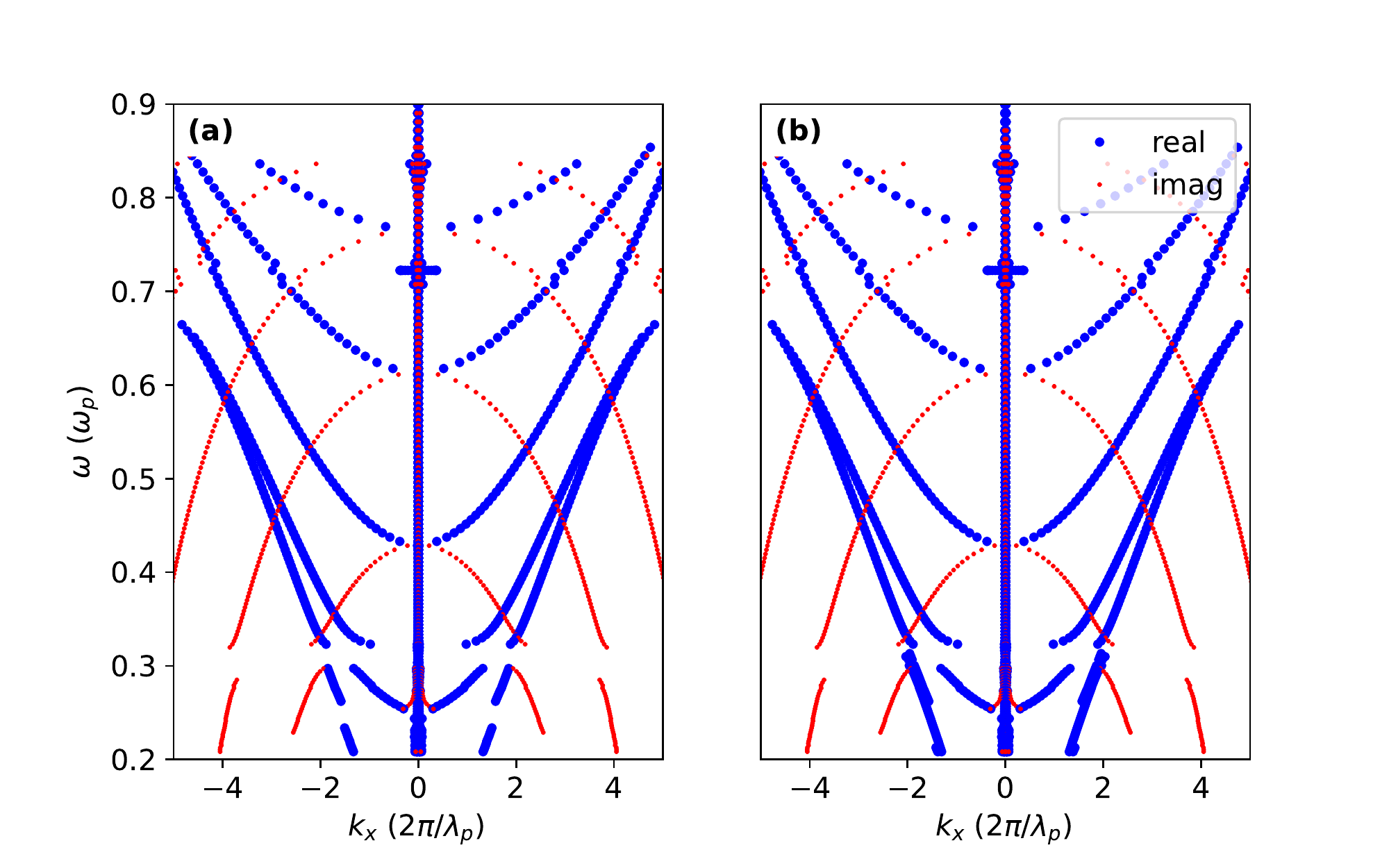}
    \caption{Band structures obtained with a) the modified Dirichlet boundary, b) with the periodic boundary condition. The blue and red dots correspond to the real and imaginary parts of the $k_x$, respectively.}
    \label{fig:band_ver}
\end{figure}

First, we validate that the use of different boundary conditions backing the PML does not affect the guided modes of the system.  In Fig. \ref{fig:mode_profiles_check}b and c, we compare the real part of the field profiles of the structure for a sample mode with the same solved $k_x$ eigenvalue but with either the Dirichlet or periodic boundary conditions. We see an excellent agreement between the modal profiles. In Fig. \ref{fig:band_ver}, we compare the bandstructure of the waveguide with the periodic and Dirichlet boundary conditions backing the PML in the $y$-direction. To generate these bandstructures, we scan a discrete set of $\omega$ values, and look for the $k_x$ eigenvalues closest to 0.  Additionally, due to spurious solutions that come from discretizing Eq. \eqref{eq:maxwell_quadratic_eigen}, a mode-filter must be implemented to produce the band-structure. Fig. \ref{fig:band_ver}a is the band structure with the Dirichlet boundary and Fig. \ref{fig:band_ver}b is the corresponding one with the periodic boundary. Again, we see that the band structures in both case are essentially identical, as expected.

\begin{figure}
    \centering
    \includegraphics[width= 5 in]{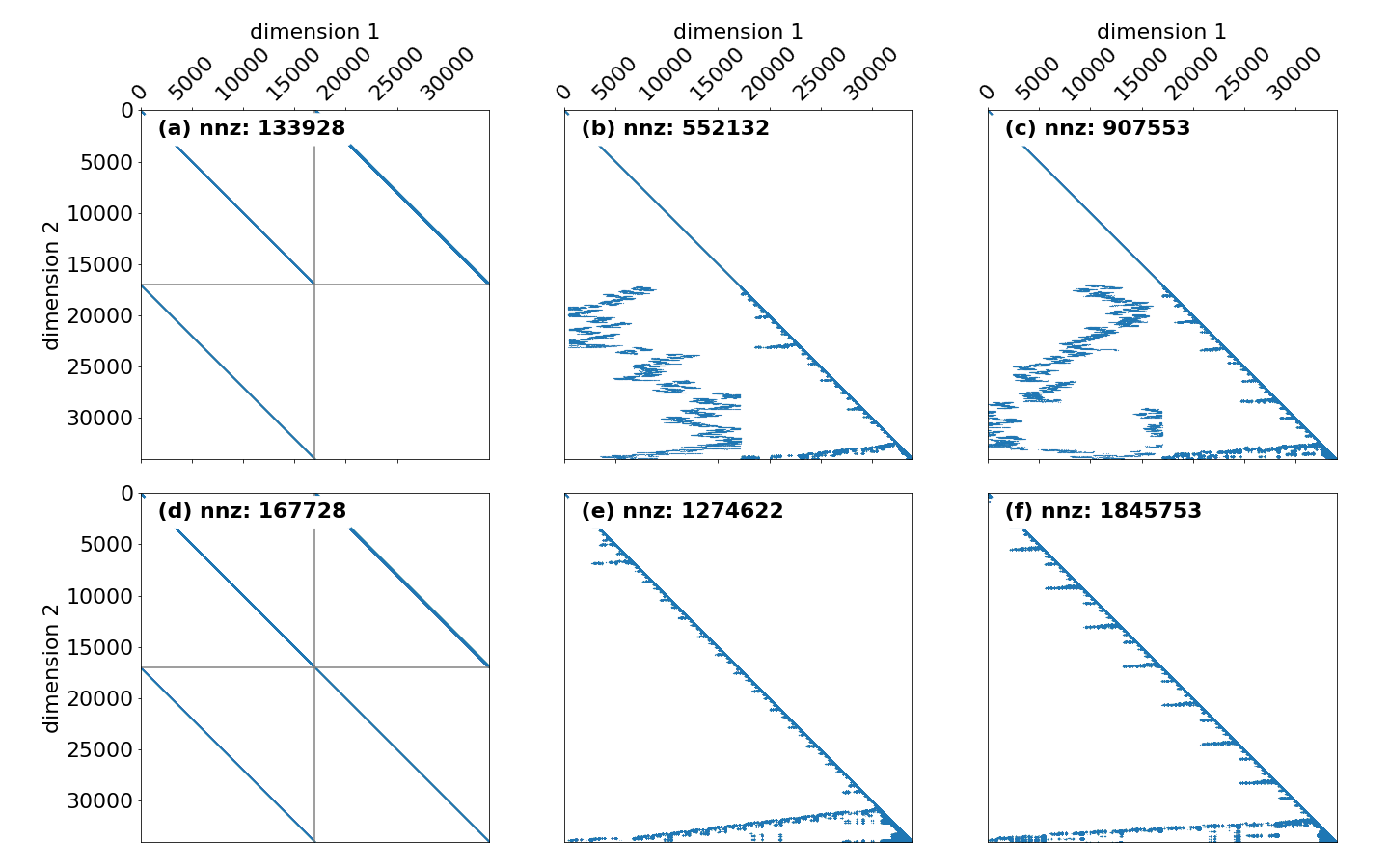}
    \caption{a) Sparsity pattern of $\hat B -\sigma \hat D$ with $\sigma = 0$ in Eq. \eqref{eq:quadratic_eigen} of the eigenvalue problem. The gray lines partition the matrix into $N_xN_y \times N_xN_y$ sub-blocks corresponding to those in Eqs. \eqref{eq:hat_b} and \eqref{eq:hat_D}. The corresponding sparsity pattern of the $\hat L$ factor with the periodic boundary condition. c) The corresponding sparsity pattern of the $\hat L$ factor for the modified Dirichlet boundary condition. In the second row, we show the analogous result with $\sigma = 1$ (in units of $\lambda_p$, the free-space wavelength at the plasma frequency).}
    \label{fig:dispersive_eigensolve}
\end{figure}

Now we compare the fill-in of factorizing $(\hat B-\sigma \hat D)$ with the periodic versus the modified Dirichlet boundary conditions. In Fig. 8a-c we consider the case where $\sigma = 0$. In Fig. \ref{fig:dispersive_eigensolve}a, we show the sparsity pattern of the matrix $\hat B$ in Eq. \eqref{eq:hat_b}.  In Fig. \ref{fig:dispersive_eigensolve}b and c, we show the fill-in of the $\hat L$ factor of $\hat B$ with the domain truncated by the periodic boundary and the modified Dirichlet boundary conditions, respectively. The fill-in patterns vary substantially between the two cases. Moreover, even though we have PML regions only along the $y$ boundaries, we still achieve a $40\%$ reduction in the number of nonzeros. In the second row of Fig. \ref{fig:dispersive_eigensolve}, we show a case with $\sigma \neq 0$. Since $\hat B$ and $\hat D$ have different sparsity patterns, the sparsity pattern shown in Fig. 8d differ from that in Fig. 8a. In comparing Fig. 8e with Fig. 8f, we again see that the use of the modified Dirichlet boundary condition leads to less fill-in. 


\section{Conclusion}

In summary, we consider direct solvers for the FDFD methods applied to Maxwell’s equations. We demonstrate that the choice of boundary condition used behind the PML can have significant effects on the numerical efficiency of the LU factorization of the system matrix. With the PML the physical properties of the solutions do not depend on the choices of the boundary conditions. But we observe that the use of a modified Dirichlet boundary condition can lead to a reduction of fill-in by up to 40\%, as compared to the commonly used periodic boundary conditions. Such a reduction translates to significant increase in computational speed and decrease in memory requirements for the direct solvers.

While in our paper we focus on the PML, our results should be applicable for other numerical absorbers, such as the adiabatic absorbers \cite{Oskooi:08}. Also, while we illustrate our observations with two-dimensional simulations, the same observations should be applicable for three-dimensional simulations as well. Our observation should also be relevant for finite-element frequency-domain solvers. 

Our result is important for enhancing the performance of direct solvers in the FDFD methods. In addition, this result should have implications for several numerical tasks in computational electromagnetics that can benefit from efficient LU decomposition of the system matrix, such as adjoint variable method\cite{Veronis2004, Lalau-Keraly2013}, efficient line searches for gradient-based photonic optimization \cite{boutami_zhao, Boutami:20, Zhao:22}, and domain decomposition methods utilizing the Schur complement approach \cite{Zhao2019,Zhao:18}.

\bibliographystyle{unsrt}  
\bibliography{references}  

\end{document}